\documentclass[12pt]{amsart}
\usepackage{amsmath,amssymb}
\newtheorem{theorem}{Theorem}
\newtheorem{proposition}[theorem]{Proposition}
\def\C{\Bbb C}
\def\D{\Bbb D}
\def\ds{\displaystyle}

\title{On a lower bound of the Kobayashi metric}

\author{Nikolai Nikolov}

\address{Institute of Mathematics and Informatics\\Bulgarian Academy
of Sciences\\ Acad. G. Bonchev 8, 1113 Sofia, Bulgaria\newline
\indent Faculty of Information Sciences\\
State University of Library Studies and Information Technologies\\
Shipchenski prohod 69A, 1574 Sofia,
Bulgaria}\email{nik@math.bas.bg}



\begin{document}

\begin{abstract} It is shown that a lower bound of the Kobayashi
metric of convex domains in $\C^n$ does not hold for non-convex domains.
\end{abstract}

\maketitle

Let $D$ be a domain in $\C^n,$ $z\in D$ and $X\in\C^n.$ Denote by
$\gamma_D$ and $\kappa_D$ the Carath\'eodory and Kobayashi metrics
of $D:$
$$\gamma_D(z;X)=\sup\{|f'(z)X|:f\in\mathcal O(D,\D)\},$$
$$\kappa_D(z;X)=\inf\{|\alpha|:\exists\varphi\in\mathcal O(\D,D)\text{
with }\varphi(0)=z,\alpha\varphi'(0)=X\},$$
where $\D$ is the unit disc. Then $\gamma_D\le\kappa_D.$
By the Lempert theorem \cite{Lem}, $\gamma_D=\kappa_D$ if $D$ is convex.

Set $d_D(z;X)$ to be the distance from $z$ to $\partial D$ in the $X$-direction,
i.e. $$d_D(z;X)=\sup\{r>0:z+\lambda X\in D\text{ if }|\lambda|<r\}$$
(possibly $d_D(z;X)=\infty$). It follows by the definition of $\kappa_D$ that
$$\kappa_D(z;X)\le\frac{1}{d_D(z;X)}.$$
On the other hand, if $D$ is convex, then
$$\gamma_D(z;X)\ge\frac{1}{2d_D(z;X)}$$
due to S. Frankel \cite[Theorem 2.2]{Fra} and I. Graham \cite[Theorem 5]{Gra}.
The following short proof of this estimate can be found in \cite[Theorem 4.1]{BP}:
$$\gamma_D(z;X)\ge\gamma_{\Pi}(z;X)=\frac{1}{2d_{\Pi}(z;X)}=\frac{1}{2d_D(z;X)},$$
where $\Pi$ is a supporting half-space of $D$ at a boundary point of
the form $z+\lambda X,$ $|\lambda|=d_D(z,X).$

It turns out that the converse to the Frankel-Graham result holds.

\begin{proposition} For a domain $D$ in $\C^n$ the following conditions
are equivalent:

\noindent (a) $D$ is convex;

\noindent (b) $\ds\gamma_D(z;X)\ge\frac{1}{2d_D(z;X)};$

\noindent (c) $\ds\liminf_{z\to a}\frac{2\kappa_D(z;z-a)-1}{||z-a||}\ge 0$
for any $a\in\partial D.$
 \end{proposition}

\noindent{\it Proof.} We have only to show that (c) implies (a). Assume the contrary.
According to \cite[Theorem 2.1.27]{Hor}, one may find a point $a\in\partial D$
and a real-valued quadratic polynomial $q$ such that $q(a)=0,$ $\nabla q(a)\neq 0,$
the set $G=\{z\in\C^n:q(z)<0\}$ is contained in $D$ near $a$ and $\partial G$ has
normal curvature $\chi<0$ in some direction $X\in T^{\Bbb R}_a\partial G.$ We may assume by
continuity that $X\not\in T^\C_a\partial G.$ Then the planar set $F=G\cap(a+\C X)$ has
smooth boundary near $a$ with curvature $\chi$ at $a.$
Let $E$ be the connected component of $F$ for which $a\in\partial E.$
Denote by $n$ the inner normal to $\partial E$ at $a.$ Using \cite[Proposition 1]{NTA},
we get the contradiction
$$0\le\limsup_{n\ni z\to a}\frac{2\kappa_D(z;z-a)-1}{||z-a||}
\le\lim_{n\ni z\to a}\frac{2\kappa_E(z;z-a)-1}{||z-a||}=\frac{\chi}{2}<0.$$

\end{document}